\documentclass{elsart}
\usepackage{amssymb}
\usepackage{amsmath}
\def\ul{\underline}
\def\ov{\overline}
\def\al{\alpha}
\def\be{\beta}

\begin{document}

\begin{frontmatter}



\title{The Qth-power algorithm in characteristic 0}


\author{Douglas A. Leonard}
\ead{leonada@auburn.edu}

\address{Department of Mathematics and Statistics\\Auburn University\\Auburn, AL 36849}

\begin{abstract}

 The {\em Qth-power algorithm}
 provides a useful canonical $P$-module presentation 
 for the integral closures of certain integral extensions of
 $P:=\mathbf{F}[x_n,\ldots,x_1]$, a polynomial ring over the finite field $\mathbf{F}:=\mathbf{Z}_q$
of $q$ elements. 
 Here it is shown how to use this for several small primes $q$ 
 to reconstruct similar integral closures over the rationals, $\mathbf{F}:=\mathbf{Q}$,
 using the {\em Chinese remainder theorem} 
 to piece together presentations in different positive characteristics,
 and the {\em extended Euclidean algorithm} 
 to reconstruct small rational fractions 
 to lift these to presentations over $\mathbf{Q}$.

\end{abstract}


\end{frontmatter}

\section{Introduction}

The {\em Qth-power algorithm}, \cite{mod},\cite{Pell},\cite{qth}, 
is designed to compute integral closures of integral extensions 
of multivariate polynomial rings in positive characteristic 
by exploiting the linearity of the Frobenius map $a\mapsto a^q$ in characteristic $q$.
Over the rationals there is no such map to exploit.

Since the integral closure over $\mathbf{Q}$ (the rationals) 
should specialize to the integral closure over $\mathbf{Z}_q$ (the finite field of $q$ elements)
for almost all primes $q$, it makes perfect sense to compute the integral closure
of the mod $q$ image for one or more primes, use the Chinese remainder theorem to reconcile
the results modulo the product of those primes, and then use the extended Euclidean algorithm
to lift those results to results in characteristic $0$.
It will be proven that the presentation over $\mathbf{Q}$ is isomorphic to a subring of the presentation of its
image over $\mathbf{Z}_q$ in general (at least when everything involved makes sense mod $q$), 
with equality for most $q$;
and that the reconciled version of the presentations over $\mathbf{Z}_q$ for various $q$ is lifted back to 
something over $\mathbf{Q}$ which has fractions with the same set of leading monomials as those of the fractions in its mod $q$ images.
So if the lifted version has an isomorphic image of the original ring inside it 
(which will not happen unless the product of the distinct primes used is large enough), 
it necessarily must be the integral closure of that original ring.  

The fact that the extended Euclidean algorithm gives essentially inverse results of the mod $q$ map 
when $q$ is sufficiently large, 
gives, as a corollary, the integral closure for large primes $q$ can be gotten from those for several smaller
primes, a useful result here, in that the Qth-power algorithm should be, by its nature, 
expected to perform significantly better for smaller primes $q$.

It should be pointed out that in both characteristic $0$ and characteristic $q>0$ the main advantage of the Qth-power
algorithm is that it takes highly structured input and produces highly structured output,
namely a strict affine P-algebra presentation for the integral closure with an induced
(as opposed to default) monomial ordering based on the weighted monomial ordering on the input.
This allows for a fairly simple determination of the existence of a better Noether normalization than the given one,
giving a smaller presentation with the same type of structure and information.
But in characteristic $0$, it gives a presentation with relatively small rational coefficients
rather than a presentation over the integers with overly large integer coefficients, and one that specializes
mod $q$ to that of the image mod $q$ for all large $q$ (and a subset of the integral closure for all smaller $q$
for which the image makes sense). 

Section 2 contains notation to describe the main algorithm, the algorithm itself, and an outline of
what will be proved to justify it. The technical definitions and other details are postponed to section 3.
Section 4 deals with the computation of canonical conductor elements based on the Jacobian.
Section 5 is a discussion of the application of the Chinese remainder theorem 
and the extended Euclidean algorithm in this context. 
And section 6 contains the theory and proofs neeeded to justify the algorithm.
There are numerous small examples throughout to help with the concepts and notation.
But there is a larger example relegated to the Appendix useful as well for comparison
to other implementations of integral closure computations. There are several other such on the author's website.
Code for Qth-power algorithm has been available for a while in {\sc Magma} and more recently
in {\sc Macaulay2} as well. The latter has documentation containing more examples.

\section{Overview of the algorithm}

Let $P^{(0)}:=\mathbf{Q}[x_n,\ldots,x_1]$ be the polynomial ring 
over the rationals in the (independent) variables $x_n,\ldots,x_1$.
Let $S^{(0)}:=P^{(0)}[y]/\langle f(y)\rangle$ for some monic polynomial $f(T)\in P[T]$,
be an integral extension. 
Suppose also that it is an affine domain of type I with a weight function
defining a weight-over-grevlex monomial ordering (as in \cite{mod}), 
though more general hypotheses may suffice.

Let $Q(S^{(0)})$ denote the field of fractions of $S^{(0)}$.
Let $\Delta^{(0)}\in P^{(0)}$ be the canonical monic conductor element
computed from Jacobian($B^{(0)}$) as described below,
so that the integral closure $C(S^{(0)},Q(S^{(0)}))$ of $S^{(0)}$ in $Q(S^{(0)})$
satisfies
\[S^{(0)}\subseteq C(S^{(0)},Q(S^{(0)}))\subseteq \frac{1}{\Delta^{(0)}}S^{(0)}\subset Q(S^{(0)}),\]
and is known to be the union of all rings lying between
$S^{(0)}$ and $\frac{1}{\Delta^{(0)}}S^{(0)}$.
[Since $S^{(0)}$ is assumed to be an integral extension of $P^{(0)}$, 
it probably makes more sense to think of it as $C(P^{(0)},Q(S^{(0)}))$, 
in that $P^{(0)}$ is a minimal subring over which this has a finite module structure. 
And once a fixed conductor element $\Delta^{(0)}\in P^{(0)}$ has been chosen,
it makes more sense to use the notation $C(P^{(0)},\frac{1}{\Delta^{(0)}}S^{(0)})$
(even though $\frac{1}{\Delta^{(0)}}S^{(0)}$ is not a ring) to emphasize that the elements
are from $\frac{1}{\Delta^{(0)}}S^{(0)}$ and are integral over $P^{(0)}$.
So that is the notation that will be used here.]

The objective here is to find a canonical ordered set of {\em monic} polynomials
(that is, with leading coefficient $1$ relative to the monomial ordering being used) 
$g_{J_0}^{(0)},\ldots,g_1^{(0)},g_0^{(0)}:=\Delta^{(0)}$ so that
the corresponding fractions $g_j^{(0)}/\Delta^{(0)}$, $0\leq j\leq J_0$
form a $P^{(0)}$-module generating set for $C(P^{(0)},\frac{1}{\Delta^{(0)}}S^{(0)})$,
then use $\ov{y}_{J_0}^{(0)},\ldots,\ov{y}_1^{(0)}$ as new variable names
for the fractions $g_{J_0}^{(0)}/\Delta^{(0)},\ldots,g_1^{(0)}/\Delta^{(0)}$
to define a ring 
\[\ov{R}^{(0)}:=\mathbf{Q}[\ov{y}_{J_0}^{(0)},\ldots,\ov{y}_1^{(0)};x_n,\ldots,x_1]\] 
with grevlex-over-weight monomial ordering induced by the weight function of $R^{(0)}$
(as in \cite{mod} but described below as well), 
and then compute the monic polynomials
$\ov{b}_k^{(0)}$, $1\leq k\leq K_0$,
defining the minimal, reduced Gr\"obner basis $\ov{B}^{(0)}$
of the ideal $\ov{I}^{(0)}$ of induced relations.

Then $\ov{S}^{(0)}:=\ov{R}^{(0)}/\ov{I}^{(0)}$
is a strict affine $P^{(0)}$-algebra presentation of
the integral closure $C(P^{(0)},\frac{1}{\Delta^{(0)}}S^{(0)})$.
While $P^{(0)}$ is an explicit subring of $\ov{S}^{(0)}$,
$S^{(0)}$ need not be;
so let $\psi^{(0)}\ :\ S^{(0)}\rightarrow\ov{S}^{(0)}$ be the inclusion map
(the identity on the $x_i$ identified in each copy of $P$ throughout, 
but not on the $y$ which is mapped
to combinations of the $\ov{y}_j^{(0)}$ and the $x_i$) so that
\[\psi^{(0)}(S^{(0)})\subseteq \ov{S}^{(0)}\subseteq \frac{1}{\Delta^{(0)}}\psi^{(0)}(S^{(0)})\subset Q(\psi^{(0)}(S^{(0)})).\]
Now consider what needs to happen for there to be an image of all this over
$\mathbf{Z}_q$ for some prime $q$, gotten by applying the mod $q$ map, $\mu_q$.

It is easy enough to define $P^{(q)}:=\mathbf{Z}_q[x_n,\ldots,x_1]$ by
identifying the variables $x_i$, and similarly to define
$R^{(q)}:=\mathbf{Z}_q[y;x_n,\ldots,x_1]$ by further identifying
the variable $y$.

If $q$ doesn't divide any denominator $\be$ of any rational fraction
$\al/\be$, $gcd(\al,\be)=1$ of any $b_k^{(0)}$, $1\leq k\leq K$,
then $B^{(q)}:=\left(b_k^{(q)}:=\mu_q(b_k^{(0)})\in R^{(q)}:\ 1\leq k\leq K\right)$
makes sense, and is still a minimal, reduced Gr\"obner basis for the ideal
$I^{(q)}$ of $R^{(q)}$ that it generates, 
though the quotient ring $S^{(q)}:=R^{(q)}/I^{(q)}$ need no longer be even a reduced ring, 
let alone an affine domain of any sort. 
[Whether it is an integral extension of $P^{(q)}$ 
can depend on how one views such things as $\mathbf{Z}_q[y_1;x_1]/\langle y_1^2\rangle$ 
in which $x_1$ doesn't appear in the defining relation.]

And if $q$ doesn't divide $LC(d^{(0)})$, 
for $d^{(0)}=LC(d^{(0)})\Delta^{(0)}$ the conductor as computed from the Jacobian of $B^{(0)}$ over $\mathbf{Z}$,
then $\Delta^{(q)}:=\mu_q(\Delta^{(0)})$
is the monic canonical conductor element that would have been computed
using the Jacobian over $\mathbf{Z}_q$.

The Qth-power algorithm is meant to work in positive characteristic
to produce a strict affine $P^{(q)}$-algebra presentation with
$P^{(q)}$-module generating set of fractions with monic numerators 
$g_{J_q}^{(q)}/\Delta^{(q)},\ldots,g_1^{(q)}/\Delta^{(q)}$,  and $g_0^{(q)}/\Delta^{(q)}:=1$,
given the variable names $\ov{y}_{J_q}^{(q)},\ldots,\ov{y}_1^{(q)}$ to define 
\[\ov{R}^{(q)}:= \mathbf{F}_q[\ov{y}_{J_q}^{(q)},\ldots,\ov{y}_1^{(q)};x_n,\ldots,x_1]\]
(having the grevlex-over-weight monomial ordering induced by that on $R^{(q)}$) 
with monic polynomials $\ov{b}_k^{(q)}$, $1\leq k\leq K_q$
forming a minimal, reduced Gr\"obner basis $\ov{B}^{(q)}$ for the
ideal $\ov{I}^{(q)}$ of induced relations, defining the presentation
$\ov{S}^{(q)}:=\ov{R}^{(q)}/\ov{I}^{(q)}$.

\newpage
The steps in the proposed characteristic $0$ algorithm based on this are then simple to understand:
\begin{alg}
\begin{enumerate}
\item Start with the finite ordered set of $($independent$)$ variables $(x_n,\ldots,x_1)$ 
defining the Noether normalization $P^{(0)}$ in characteristic $0$, 
the (dependent) variable name $y$ used to define the ring $R^{(0)}$,
and the finite ordered set of monic relations $(b_1,\ldots,b_K)$ 
forming a minimal, reduced Gr\"obner basis for the ideal of relations $I^{(0)}$
for a presentation of the input quotient ring $S^{(0)}=R^{(0)}/I^{(0)}$.

\item Compute a canonical conductor element $\Delta^{(0)}\in P^{(0)}$ for $S^{(0)}$
from the Jacobian .
 
\item For successive primes, $q_l$, test that the $(mod\ q_l)$ map, $\mu_{q_l}$, is defined 
$($that is that $q_l$ doesn't divide $\be$ 
for any rational coefficient $\al/\be$, $gcd(\al,\be)=1$, 
in any of the basis relations $b_k$ or of $\Delta^{(0)})$
and that $S^{(q_l)}$ is still an integral extension of the image $P^{(q_l)}$. 

\item Compute a canonical conductor element $\Delta^{(q_l)}\in P^{(q_l)}$ for $S^{(q_l)}$,
skipping $q_l$ if it is one of the $($finite number of$)$ primes for which 
$\Delta^{(q_l)}\neq \mu_{q_l}(\Delta^{(0)})$. 

\item Use the Qth-power algorithm in characteristic $q_l$ 
$($as a black box for the purposes of this paper$)$ 
to compute a canonical ordered set of $($numerator$)$ polynomials $(g_j^{(q_l)})$, $0\leq j\leq J_{q_l}$ 
$($with the common denominator polynomial being $g_0^{(q_l)}:=\Delta^ {(q_l)})$
for the fractions forming a $P^{(q_l)}$-module generating set for a 
strict affine $P^{(q_l)}$-algebra presentation of the integral closure $C(P^{(q_l)},\frac{1}{\Delta^{(q_l)}}S^{(q_l)})$
with $\left(\ov{b}_k^{(q_l)}\right)_{l\in L}$, $1\leq k\leq K_{q_l}$, the minimal, reduced Gr\"obner basis
for the ideal of relations relative to the induced grevlex-over-weight monomial ordering.

\item If $(q_l)$, $l\in L$ is a sequence of distinct primes for which presentations
$\ov{S}^{(q_l)}$ have been computed, and 
\begin{itemize}
\item $J_L:=J_{q_l}$ is independent of $l\in L$; 
\item $LM(g_j^{(q_l)})$ is independent of $l\in L$;
\item $K_L:=K_{q_l}$ is independent of $l\in L$;
\item $LM(\ov{b}_k^{(q_l)})$ is independent of $l\in L$;
\end{itemize}
then use the Chinese remainder theorem 
on the canonical ordered sets $\left(g_j^{(q_l)}\right)_{l\in L}$, $1\leq j\leq J_L$, 
and also on the sequences $\left(\ov{b}_k^{(q_l)}\right)_{l\in L}$, $1\leq k\leq K_L$,
to get similar canonical ordered sets
 $(g_j^{(N_L)})$, $1\leq j\leq J_L$,
and  sequences $\left(\ov{b}_k^{(N_L)}\right)$, $1\leq k\leq K_L$,
for  $N_L:=\prod \{q_l:\ l\in L\}$
monic with the same sets of leading monomials.

\item Then use the extended Euclidean algorithm to lift the coefficients $c(mod\ N_L)$ 
to small fractions $\al/\be\in\mathbf{Q}$ with $\al^2+\be^2$ minimal,
to get a canonical ordered set of polynomials
$(g_j^{(0,N_L)})$, $1\leq j\leq J_L$, and sequence $(\ov{b}_k^{(0,N_L)})$, $1\leq k\leq K_L$, 
over the rationals with the same sets of leading monomials to describe 
a possible integral closure in characteristic $0$.

\item Stop when 
$(\ov{b}_k^{(0,N_L)})$, $1\leq k\leq K_L$, 
is a minimal, reduced Gr\"obner basis for the ideal $\ov{I}^{(0,N_L)}$ it generates 
and the image $\psi^{(N_L)}(I^{(0)})$ of the original ideal $I^{(0)}$  is contained in $\ov{I}^{(0,N_L)}$,
both necessary conditions. 
\end{enumerate}
\end{alg}

To explain why this works, the important steps are to show that whenever everything involved makes sense
\begin{enumerate}
   \item \[ \mu_q\left(C(P^{(0)},\frac{1}{\Delta^{(0)}}S^{(0)})\right)\subseteq C(P^{(q)},\frac{1}{\Delta^{(q)}}S^{(q)}) \]
   \item $\ov{S}^{(0,N_L)}$ is isomorphic to $\ov{S}^{(0)}$ if $I^{(0)}$ is isomorphic to a subideal of $\ov{I}^{(0,N_L)}$ and $\ov{B}^{(0,N_L)}:=(\ov{b}_k^{(0,N_L)})$, $1\leq k\leq K_L$,  is still a minimal, reduced Gr\"obner basis for $\ov{I}^{(0,N_L)}$.
\end{enumerate}
These are proven as Lemma 9 and Theorem 15 below.
It should be noted that (regardless of characteristic) any $P$-module between $S$ and $\frac{1}{\Delta}S$
has a canonical ordered set of polynomials $(g_j)_j$ (as defined in the next section) 
with $\langle LM(g_j)_j\rangle$ a measure of the size of the $P$-module.
That is, were the de Jong algorithm implemented relative to a fixed Noether normalization $P$ and a canonical
conductor element $\Delta\in P$, 
then the sequence of nested rings produced 
would have a sequence of canonical ordered sets of polynomials $(g_j)_j$ with  $\langle LM(g_j)_j\rangle$
nested and getting larger. 
The reverse is true of the Qth-power algorithm 
in that the $P$-modules produced have canonical ordered sets of polynomials $(g_j)_j$ 
with  $\langle LM(g_j)_j\rangle$ nested and getting smaller. 
Both approaches meet in the middle with a $P$-module that is a ring that must be the integral closure sought.

Also, Jacobian($B^{(0)}$), over the integers, $\mathbf{Z}$, is used below to define canonical
conductor elements $\Delta^{(0)}\in P^{(0)}$ and $\Delta^{(q)}\in P^{(q)}$ for all primes $q$ at the same time.  
This is discussed in its own section.
The use of the Chinese remainder theorem and the extended Euclidean algorithm,
while discussed below as applied in this context, are assumed to be elementary.
The proof that the Qth-power algorithm works in positive characteristic
was dealt with in the author's previous papers cited in the introduction, though
certain parts of it are discussed below. [It should be noted however,
that it would not take too much work to put other implementations of other
integral closure algorithms in a form that would also work here, though
at present, few if any give a similar canonical result in characteristic $0$
that directly specializes to the result they give in positive characteristic.
It makes mathematical sense to rewrite them to reflect this connection
between integral closures in characteristic $0$ and positive characteristic $q$.] 

\newpage
\section{Definitions and other details}

The following material describes the {\em structure} that is used to describe
integral closures of integral extensions of a given Noether normalization $P$,
and explaining the mindset of the Qth-power algorithm approach to same.
The idea is to have an integral extension $S:=R/I$ of $P:=\mathbf{F}[x_n,\ldots,x_1]$, 
with $R:=\mathbf{F}[y;x_n,\ldots,x_1]$, 
$I:=\langle f(y)\rangle$ the ideal of relations, 
and $P$ a Noether normalization of $S$,
and with $S$ having a weight function induced on it by $P$, \cite{mod}.
Then its integral closure $C(P,S)$ 
should have a presentation $\ov{S}:=\ov{R}/\ov{I}$
with an induced weight function. 
That is, for $\ov{R}:=\mathbf{F}[\ov{y}_J,\ldots,\ov{y}_1,x_n,\ldots,x_1]$ with
$\ov{y}_j$ a name for the (non-trivial) $P$-module generator $g_j/g_0$ having
$wt(\ov{y}_j):=wt(g_j)-wt(g_0)$ as its induced weight.

Moreover a presentation $\ov{S}$ of the integral closure $C(S,Q(S))$ should have a nice structure as an affine $P$-algebra. 

\begin{defn} A {\em strict} affine $P$-algebra presentation $\ov{S}:=\ov{R}/\ov{I}$
with $\ov{R}:=\mathbf{F}[\ul{\ov{y}},\ul{x}]$, 
and $\ov{I}$ the ideal of induced relations
is one with a minimal, reduced Gr\"obner basis $\ov{B}$ for $\ov{I}$ 
consisting of $P$-quadratic relations of the form 
$\ov{y}_i\ov{y}_j-\sum_kc_{i,j,k}\ov{y}_k$, $c_{i,j,k}\in P$, 
describing the $P$-algebra multiplication
with possibly some monic $P$-linear relations of the form 
$\sum_ka_k\ov{y}_k$, $a_k\in P$,  
if the $P$-module generators, $\ov{y}_k$, are not independent over $P$, \cite{mod}.
\end{defn}

This is ensured by the grevlex-over-weight monomial ordering, 
in that all products of total degree 2 in the $\ov{y}$'s 
are reduced to $P$-linear combinations of the $P$-module generators 
with total degree less than two, 
and all $P$-syzygies only involve monomials of total degree less than 2 in them.

\begin{defn}
For $P:=\mathbf{F}[x_n,\ldots,x_1]$ a polynomial ring, $R:=\mathbf{F}[y;x_n,\ldots,x_1]$, $I$ an ideal of $R$ 
such that the quotient ring $S:=R/I$ is an integral extension of $P$, an ordered set $(g_j\in R: 0\leq j\leq J)$
of polynomials is said to be {\em canonical} for some submodule $\frac{1}{\Delta}T\subseteq \frac{1}{\Delta}S$ iff
\begin{enumerate}
   \item each $g_j$ is monic $($has leading coefficient $1$ relative to the monomial ordering being used$)$;
   \item $g_0=\Delta\in P$ is a conductor element for $T$;
   \item $\left(g_j/\Delta\ :\ 0\leq j\leq J\right)$ is a $P$-module generating set for $\frac{1}{\Delta} T$;
   \item $(g_j\ :\ 0\leq j\leq J)$ is {\em interreduced}, meaning that no monomial of any $g_{j_1}$
         is of the form $\ul{x}^{\ul{\alpha}}LM(g_{j_2})$ for any $j_2\neq j_1$ and  $\ul{x}^{\ul{\alpha}}\in P$.
\end{enumerate}
\end{defn}

Making sure that one produces a canonical strict affine $P$-algebra presentation, independent of characteristic,
is crucial in being able to reconstruct a canonical presentation in characteristic 0 from canonical one in positive
characteristic.

Moreover the {\em induced grevlex-over-weight monomial ordering}
is based on a {\em weight function} on the input extended to the output: 

\begin{defn} A {\em weight function} $wt\ :\ R\backslash I\rightarrow \mathbf{N}^n$ $(n$ the number of independent variables$)$
is a function satisfying:
\begin{enumerate}
\item $wt(c)=\ul{0}\mbox{ for all }c\in\mathbf{F}\backslash \{ 0\}$;
\item $wt(gh)=wt(g)+wt(h)$ for all $g,h$;
\item if $wt(g)=wt(h)$, then for some unique $c\in \mathbf{F}$ either $wt(g-ch)\prec_{lex} wt(g)$ or $g-ch\in I$.
\end{enumerate}
\end{defn}

If $R=P=\mathbf{F}[x_n,\ldots,x_1]$, the columns of any non-singular matrix $M_P$ 
defining a (global) monomial ordering on $P$ 
also define a weight function on (the non-zero elements of) $P$.

\begin{exmp} $\mathbf{F}[y,x]/\langle y^3-x^5-yx\rangle$ has a weight function 
with $wt(y):=5$ and $wt(x):=3$ $($corresponding to pole orders
of these as rational functions at the only point where either has a pole$)$. 
The example
$\mathbf{F}[y,x]/\langle y^3-x^3y-x\rangle$ does not. $[$An attempt
at weights would have $wt(y^3)=wt(x^3y)$, so $3wt(y)=wt(y)+3wt(x)>wt(x)>0$
with no field element $c$ with $wt(y^2-cx^3)$ a smaller weight
because $y^2-x^3=x/y$ would have weight less than $0$.$]$
The example $\mathbf{F}[x,y,z]/\langle x^6+x^3z-y^2z\rangle$ has
a weight function with $wt(x):=(3,2)$, $wt(y):=(6,6)$, $wt(z):=(6,0)$, 
an extension of the grevlex monomial order with $wt(y)=(1,1)$ and $wt(z)=(1,0)$.
\end{exmp}

Weight functions (relative to a given ideal $I$) have the important property 
that $wt(g)=wt(NF(g,I))$, as otherwise their difference (an element of $I$) would have a defined weight. 
This, in turn, implies that all standard monomials have different weights.
Integral extensions with weight functions have at least this much more structure than those that don't. 
The integral extensions such as those considered here have a weight function.

One can extend a weight function on $S$ naturally 
to (the non-zero elements of) $Q(S)$ 
(the {\em field of fractions} of $S$) 
by $wt_{Q(S)}(g/h):=wt_S(g)-wt_S(h)$ if one allows values in $\mathbf{Z}^n$.
But for $g/h\in C(S,Q(S))$,  
$wt(g/h)$ will not have negative entries and will represent an induced weight on $g/h$.



There are weight-over-grevlex and grevlex-over-weight monomial orderings
defined by the non-singular matrices gotten by replacing the top $n$ or bottom
$n$ rows of a grevlex monomial ordering matrix by the weight matrix, respectively
The former emphasizes the property that $wt(LT(f))=wt(LT(f-LT(f)))$,
whereas the latter emphsizes the desired strict affine $P$-algebra presentation.

The conductor element $\Delta$ and the numerators $g_j$ produced 
are all assumed to be {\em monic}. 
What will be computed are pairs of finite canonical ordered sets of polynomials 
and finite sequences of relations, 
with maps between such pairs. 
The induced weights are kept track of as well, 
given that they define the monomial orderings involved.

What is necessary to know about the {\em Qth-power algorithm} is that it treats 
the input ring $S$ as a $P$-module with a natural induced monomial ordering, 
computes a {\em conductor element} $\Delta\in P$,
starts with a dual module such as the default $M_0:=\frac{1}{\Delta}S$ 
and computes a nested sequence of $P$-modules
\[M_0\supset\cdots\supset M_L=M_{L+1}=C(P,\frac{1}{\Delta}S)\]
by the simple definition
\[M_{i+1}:=\left\{\frac{g}{\Delta}\in M_i:\ \frac{NormalForm(g^q,I)}{\Delta^q}\in M_i\right\}.\]
Necessarily each $M_i$ (and hence the integral closure itself) 
is a $P$-module with a natural induced monomial ordering.
Moreover it naturally produces a strict affine $P$-algebra presentation 
$\ov{R}/\ov{I}$ relative to a canonical ordered set $(g_j)_j$ of polynomials. 
[This approach works theoretically for any characteristic 
and any integer power at least $2$,
but is only linear when $q$ is (a power of) the characteristic.]

Consider the reasonably generic example:
\begin{exmp} Let $P^{(0)}:=\mathbf{Q}[x]$, $R^{(0)}:=\mathbf{Q}[y;x]$, and $I^{(0)}:=\langle B^{(0)}\rangle$, with
\[B^{(0)}:=\{(y^2-3/4y-15/17x)^3-9yx^4(y^2-3/4y-15/17x)-27x^{11}\},\]
with weight function $(11,6)$ 
$($meaning that the monomials $y^6$ and $x^{11}$ are the only two in the defining relation of largest weight, $66$$)$
defining an integral extension $S^{(0)}=R^{(0)}/I^{(0)}$. 
$[$This was originally constructed to have at least one non-trivial integral element,
namely $(y^2-3/4y-15/17x)/x^2$ and several moderate-sized rational coefficients
to be reconstructed, in order to test this extension of the Qth-power algorithm to characteristic $0$.$]$

This is worked out in detail in the Appendix, not only by the methods being described here
and implemented in the author's {\em QthPower} package in {\sc Macaulay2} \cite{JSAG},
but also using the other existing applicable implementations of integral closure and/or normalization
algorithms, {\em normal} in {\sc Singular}, {\em integralClosure} in {\sc Macaulay2},
and both {\em Normalisation} and {\em IntegralClosure} in {\sc Magma}.
\newpage

The output from the {\em qthIntegralClosure} function in the {\em QthPower} package consists of:
\begin{enumerate}
\item an ordered set of ``numerators'' for a $P^{(0)}$-module generating set $($here a basis$)$:
\[g_0^{(0)}/x^4=\delta^{(0)}=\Delta^{(0)}/x^4=x^5,\] 
\[g_1^{(0)}/x^4=y^2x^3-\frac{3}{4}yx^3-\frac{15}{17}x^4,\] 
\[g_2^{(0)}/x^4=yx^5,\]
\[g_3^{(0)}/x^4=y^4x-\frac{3}{2}y^3x-\frac{30}{17}y^2x^2+\frac{9}{16}y^2x+\frac{45}{34}yx^2+\frac{225}{289}x^3,\] 
\[g_4^{(0)}/x^4=y^3x^3-\frac{15}{17}yx^4-\frac{9}{16}yx^3-\frac{45}{68}x^4,\]
\[g_5^{(0)}/x^4=y^5-\frac{9}{4}y^4-\frac{30}{17}y^3x+\frac{27}{16}y^3+\frac{45}{17}y^2x+\frac{225}{289}yx^2-\frac{27}{64}y^2-\frac{135}{136}yx-\frac{675}{1156}x^2\]
with $g_0^{(0)}/x^4=\delta^{(0)}=x^5\in P^{(0)}$ also being the $($reduced$)$ conductor element and common denominator of the fractions, 
$($not the actual $\Delta^{(0)}=x^9$, accessible only through the {\em qthConductor} function$)$;
\item a minimal, reduced Gr\"obner basis of induced relations;
\item an induced weight matrix $($here a weight matrix$)$ $(25,21,20,11,10,6)$
defining a grevlex-over-weight monomial ordering.
\end{enumerate}

A typical {\em induced} relation in the Gr\"obner basis, such as
\[ \ov{y}_4\ov{y}_2-3/2\ov{y}_4-\ov{y}_3x_1^2-15/17\ov{y}_1x_1+9/8\ov{y}_1,\]
$($produced by finding the normal form of $\ov{y}_4\ov{y}_2$ relative to the input ideal $I^{(0)}$$)$
corresponds to some $P$-algebra multiplication rule, in this case
\[ f_{21}*f_{11}=3/2f_{21}+f_6^2f_{20}+\left(15/17f_6-9/8\right)f_{10}\]
$($using notation $f_{25}:=\ov{y}_5$, $f_{21}:=\ov{y}_4$, $f_{20}:=\ov{y}_3$, $f_{11}:=\ov{y}_2$, 
$f_{10}:=\ov{y}_1$, $f_6:=x_1$ to reflect the weights$)$.
Note that the weight of the left side, $f_{21}f_{11}$ and the weight of its
NormalForm, the right-hand side, are both necessarily the same $($in this case, $32$$)$.
\end{exmp}

\newpage
\section{Computing a canonical conductor element}

Standard methods to compute a {\em conductor element} $\Delta\in S$ 
(meaning an element for which $C(S,Q(S))\subseteq \frac{1}{\Delta} S$) 
use determinants of $n\times n$ minors of a Jacobian.
This can be easily done by column-reducing the Jacobian matrix $Jacobian(B)$ of $B$;
and this computation can be done over $R$ instead of $S$ by appending columns
one for each basis element of $I$ and each row of $Jacobian(B)$.
It is then possible to consider those entries $C_{i,j}\in P$ of the column-reduced form, $C$, 
for which $C_{k,l}=0$ for $k>i$ and $l\leq j$.
An appropriate monomial ordering must be chosen relative to which this is done, 
so that the elements $C_{i,j}\in P$ considered will 
correspond to diagonal entries in $n\times n$ minors 
whose determinants necessarily produce conductor elements, 
greatest common divisors of those in the same row can be used,
and a scaled product of those gcds over all rows can be used to give a canonical conductor element $\Delta\in P$.
For this purpose, any block ordering treating the dependent variables any way 
but using the given monomial ordering described by $M_P$ 
on the lowest block consisting of the (independent) variables in $P$, will suffice.
[Note that when computing $\Delta^{(0)}\in P^{(0)}$ over $\mathbf{Q}$,
it is possible to do this over the integers, $\mathbf{Z}$, instead (if denominators are cleared first)
in order to see in one computation for what (finite set of) primes $q$ it might be that $\Delta^{(q)}\neq \Delta^{(0)}(mod\ q)$,
by seeing what primes occur anywhere in the column reduction $C$.]

The method {\em qthConductor} exported from the author's {\em QthPower} package, \cite{JSAG}, in {\sc Macaulay2}  
can be used to compute such a canonical conductor element, by letting {\sc Macaulay2} do the column-reduction,
then using a simple loop to compute the product of the gcds described. 
This computation is not a point of this paper, other than to insure that
there is a canonical conductor element that can be computed,
that it is an element of the given {\em Noether normalization}, $P$,
and that the computation in positive characteristic mirrors the computation
in characteristic $0$.

Consider the following instructive example,
meant originally to test minimality and form of presentation,
but, as a byproduct, was used to catch bugs in various implementations as well.
\begin{exmp}
Let $P^{(0)}:=\mathbf{Q}[x]$, $f^{(0)}(T):=T^8-T^2x^3+2Tx^6-x^9\in P^{(0)}[T]$,
$R^{(0)}:=\mathbf{Q}[y,x]$ with a monomial ordering based on 
$W_S:=(9,8)$,
and $S^{(0)}:=R^{(0)}/\langle f^{(0)}(y)\rangle$.
The $($extended$)$ {\em Jacobian matrix}
$(\partial f^{(0)}(y)/\partial y,\ \partial f^{(0)}(y)/\partial x,\ f^{(0)}(y))$
column-reduces over $\mathbf{Z} $ 
with {\em lex}, $x\prec y$, monomial ordering to
\[(120x^{24},\ 15x^{26},\ 5x^{27},\ x^{41}+6x^{26},\]
\[ 6yx^5-x^{38}-24x^{23}-6x^8,\ 2yx^6-2x^{39}-48x^{24}-2x^9,\ yx^{26}+2x^{29},\]
\[ 3y^2x^2-2x^{38}-48x^{23}-3x^8,\ y^2x^6-2x^{27}-x^{12},\  8y^7-2yx^3+2x^6,\]
\[ y^7x^2+yx^{20}+2yx^5+x^{38}-10x^{23}-2x^8,\ y^8-y^2x^3+2x^{39}+48x^{24}+x^9).\]
From this it is easy to extract a {\em conductor element} $\Delta^{(0)}:=x^{24}\in P^{(0)}$. 
It is also easy to extract $\Delta^{(q)}:=x^{24}\in P^{(q)}$ for $q\neq 2,3,5$,
$\Delta^{(2)}:=x^{26}$, $\Delta^{(3)}:=x^{27}$, and $\Delta^{(5)}:=x^{26}(x^3+1)^5$
for similar problems with $\mathbf{F}=\mathbf{Z}_q$, $q$ a prime.

There are rational functions $(f_0:=1,f_4,f_5,f_9,f_{10},f_{14},f_{15},f_{19})$
$($with the subscripts corresponding to the weights$)$
forming a $P^{(0)}$-module basis for the integral closure, $\ov{S}^{(0)}$.
Then the presentation of $\ov{S}^{(q)}$ can be gotten by reading $\ov{S}^{(0)}$ modulo $q$
for all primes $q\neq 3,5$.
Curiously, the smallest conductor element that could be used is $\delta^{(0)}=\delta^{(q)}=x^{13}$ 
for all primes except $\delta^{(5)}=x^{13}(x^3+1)^2$.
It is tempting to conjecture that $\Delta^{(q)}=\Delta^{(0)}(mod\ q)$ implies that
$\ov{S}^{(q)}=\ov{S}^{(0)}(mod\ q)$.
It is clearly not true that
 $\delta^{(q)}=\delta^{(0)}(mod\ q)$ implies that
$\ov{S}^{(q)}=\ov{S}^{(0)}(mod\ q)$ from $q=3$ in this example;
and it is clearly not true that
$\ov{S}^{(q)}=\ov{S}^{(0)}(mod\ q)$ implies that $\Delta^{(q)}=\Delta^{(0)}(mod\ q)$
from $q=2$.
\end{exmp}

Since it is computationally easy to avoid all the (finitely many) primes $q$ 
for which $\Delta^{(q)}\neq \Delta^{(0)}(mod\ q)$ 
(necessarily divisors of some coefficient in the computation over $\mathbf{Z}$),
it is possible to simplify subsequent computations by so doing.
\newpage
The exportable {\tt QthPower} code in {\sc Macaulay2} for this is \cite{JSAG}: 
\begin{alg}
\begin{verbatim}
qthConductor = method(TypicalValue => RingElement);
qthConductor(Ideal,ZZ) := (I,depno) -> (
    R := ring I;
    RP := (coefficientRing R)[gens R,MonomialOrder=>
         {Position=>Up,{depno,#gens R-depno}}];
    IP := sub(I,RP);
    GP := gens gb (transpose jacobian IP|
         matrix{{gens IP}}**identity(RP^(numColumns(jacobian IP))));
    depvars:=take(gens RP,depno);
    indvars:=take(gens RP,depno-#gens RP);
    qthconductor := 1;
    rowconductor := 0;
    j := numColumns(GP)-1;
    i := numRows(GP)-1;
    while i >= 0 and j >= 0 do(
       while i>=0 and j>=0 and (GP_(i,j) == 0 
       or (logpoly(GP_(i,j),depvars,indvars))#1 != 1) do(
          j = j-1;
       );
       rowconductor = 0;
       while i >= 0 and j >= 0 and GP_(i,j) != 0 
       and (logpoly(GP_(i,j),depvars,indvars))#1 == 1 do(
         rowconductor=gcd(rowconductor,GP_(i,j));
         j = j-1;
       );
       if rowconductor != 0 then qthconductor = qthconductor*rowconductor;
       i = i-1; 
    );
    s:=sub(qthconductor,R);
    s/leadCoefficient(s)
);     
\end{verbatim}
\end{alg}

\newpage
\section{Chinese remainder theorem and extended Euclidean algorithm}

As stated above, the approach given in this paper is, in some sense, an elementary one,
in that it exploits commonly known information from the Chinese remainder theorem,
and intermediate information given by the extended Euclidean algorithm.
Even a good reference such as \cite{GandG} doesn't necessarily use the
extended Euclidean algorithm in this exact way.
And, as a warning, this approach is extremely tricky in the sense that the
maps are not homomorphisms of the whole rings involved, but do extend naturally
to ring homomorphisms when defined correctly on the finite ordered sets of objects
used to describe those rings.
That is, a presentation of the integral closure of an affine $P$-algebra can be 
described in terms of a {\em finite} ordered set of fractions 
and a {\em finite} ordered set of relations.
It is then possible to construct those two {\em finite} ordered sets,
define mappings,
and then extend those naturally to definitions of ring homomorphisms.

Note especially that in implementing this approach, 
care must be taken to assure that the integral closure algorithm produces 
the same canonical result for each good prime $q$. 
That is, for most primes $q$, the integral closure over $\mathbf{Z}_q$
should look exactly like that of the integral closure over $\mathbf{Q}$,
but with coefficients reduced mod $q$.

For each presentation of $S:=R/I$ 
and presentation of its integral closure $\ov{S}:=\ov{R}/\ov{I}$, 
there is a map $\psi:R\rightarrow \ov{R}$, 
necessarily with $\psi(I)\subseteq \ov{I}$,
so that $\psi$ can be viewed as an {\em inclusion map}
$\ov{\psi}:S\rightarrow \ov{S}$.

The {\em extended Euclidean algorithm} can be used to move between
$\frac{\al}{\be}\in\mathbf{Q}$, with $gcd(\al,\be)=1,\ \be>0$ and representatives $c\in\mathbf{Z}_N$.

\begin{defn} The {\em rational fraction reconstruction map} $($see, for instance, \cite{GandG}$)$ is
\[ E_N(c):=\frac{\al}{\be},\ c\be \equiv \al(mod\ N),\ \al^2+\be^2\mbox{ minimum },\ \be>0\mbox{ minimum}.\]
The {\em mod N map} is
\[\mu_N\left(\frac{\al}{\be}\right):=c,\ c\be\equiv \al(mod\ N),\ |c|\mbox{ minimum }.\]
\end{defn}

These are almost inverse operations in the sense that
for $-N/2<c<N/2$, $\mu_N\circ E_N(c)=c$; while, for $\al^2+\be^2<N$,
$E_N\circ \mu_N(\frac{\al}{\be})=\frac{\al}{\be}$.

Both maps naturally extend to polynomials, by applying them to coefficients
and mapping variables to corresponding variables; so we shall abuse
notation and use the same function names when applying them to polynomials.

\begin{defn}
Similarly the {\em Chinese remainder map} 
standardly takes ordered sets of remainders $(a_l)_{l\in L}$ 
and ordered sets of respective moduli $(q_l)_{l\in L}$, 
and produces $a(mod\ N_L)$ for $N_L:=\prod\{q_l:\ l\in L\})$ 
such that $a\equiv a_l(mod\ q_l)$ for all $l\in L$ when the moduli are all relatively prime, 
as they will necessarily be here when the $q_l$ are distinct primes.
\end{defn}

We shall call this map {\em CRT} regardless of the number of inputs, and regardless
of whether we are applying it to integers or extending it to polynomials
by applying it to the coefficients.

\begin{lem} $\mu_q\left(C(P^{(0)},\frac{1}{\Delta^{(0)}}S^{(0)})\right)\subseteq C(P^{(q)},\frac{1}{\Delta^{(q)}}S^{(q)})$ 
for all primes $q$ for which $\mu_q\left(C(P^{(0)},\frac{1}{\Delta^{(0)}}S^{(0)})\right)$ makes sense.
\end{lem}

\begin{pf}
For each fraction $g_j^{(0)}/\Delta^{(0)}$ in the desired integral closure over $\mathbf{Q}$,
let $f_j^{(0)}(T)\in P^{(0)}[T]$ be a monic polynomial satisfied by it.
If both $\mu_q\left(g_j^{(0)}\right)/\mu_q\left(\Delta^{(0)}\right)$ 
and $f_j^{(q)}(T):=\mu_q\left(f_j^{(0)}[T]\right)\in P^{(q)}[T]$ are defined
(meaning the the prime $q$ doesn't divide the denominator $\be$ of any rational fraction
$\al/\be$, $gcd(\al,\be)=1$ occurring in either $g_j^{(0)}$ or $f_j^{(0)}[T]$),
then $\mu_q\left(g_j^{(0)}\right)/\mu_q\left(\Delta^{(0)}\right)\in \frac{1}{\Delta^{(q)}}S^{(q)}$ 
satisfies the monic polynomial $f_j^{(q)}(T)\in P^{(q)}[T]$.
\end{pf}

\begin{defn}
A prime $q$ is a {\em good} prime iff 
\[\mu_q\left(C(P^{(0)},\frac{1}{\Delta^{(0)}}S^{(0)})\right)=C(P^{(q)},\frac{1}{\Delta^{(q)}}S^{(q)}).\] 
\end{defn}

\begin{cor}
If $(g_j^{(0)} : 0\leq j\leq J)$ is canonical for $C(P^{(0)},\frac{1}{\Delta^{(0)}}S^{(0)})$
and the mod $q$ map makes sense, then 
$(g_j^{(q)} : 0\leq j\leq J)$ is canonical for $C(P^{(q)},\frac{1}{\Delta^{(q)}}S^{(q)})$
if $q$ is a good prime $($and only for some subring if it is not a good prime$)$.
\end{cor}

Clearly if $q$ divides any denominator of any rational coefficient $\al/\be$ 
of any term of any $b_k^{(0)}$, it is not good.
If $q$ divides any numerator of any rational coefficient $\al/\be$ 
of any term of any $b_k^{(0)}$, it is may not be good, 
especially if the extension $mod\ q$ is no longer really an extension.
And if $\Delta^{(q)}\neq \mu_q\left(\Delta^{(0)}\right)$, $q$ may not be good.
So computationally one can try to avoid such primes that are not good or may not be good 
(since these form a finite predictable set of primes). 

\begin{exmp} 
Consider the example with $B^{(0)}:=\{ y_1^2+13/22(x_1^9+x_1^7+x_1^5)\}$,
for which we should expect $\ov{y}_1^2=y_1/x_1^2$ and $\ov{B}^{(0)}=\{\ov{y}_1^2+13/22(x_1^5+x_1^3+x_1)\}$.
The primes $q=2,11$ are clearly bad since the $mod\ q$ map, $\mu_q$, makes no sense;
but $q=13$ is also bad in the sense that $B^{(13)}=\{y_1^2\}$ really doesn't define an integral extension
of $P^{(13)}=\mathbf{Z}_{13}[x_1]$.
Column-reducing the Jacobian $(22y_1^2+13(x_1^9+x_1^7+x_1^5),44y_1,13(9x_1^8+7x_1^6+5x_1^4))$ 
over $\mathbf{Z}$ with lex ordering gives 
$(780x_1^4,156x_1^5,52x_1^6+260x_1^4,13x_1^8+39x_1^6+65x_1^4,44y_1,22y_1^2+26x_1^7+52x_1^5)$.
So the canonical conductor elements are $\Delta^{(0)}=\Delta^{(q)}=x_1^4$ for all other $q$ except
$\Delta^{(3)}=x_1^6-x_1^4$ and $\Delta^{(5)}=x_1^5$. $q=5$ happens to be a good prime in this example,
but for $q=3$, $\ov{y}_1=y_1/(x_1^2(x_1^2-1))$, meaning there is a larger than expected integral closure $\ov{S}^{(3)}$.
Avoiding the primes $2,3,5,11,13$ $($whether or not they are not good$)$, using $7,17,19$ $($which should be good$)$ is enough to reconstruct $\ov{S}^{(0)}$, since
$7\cdot 17\cdot 19=2261>22^2+13^2=653$.
\end{exmp}

The {\em Euclidean algorithm}, applied to $N_L:=r_{-1}$ and any $r_0>0$, 
produces sequences $(r_i)$ and $(Q_i)$ such that $r_{i-2}=Q_ir_{i-1}+r_i$ with
$0\leq r_i<r_{i-1}$, and $r_n=0$.
Part of the {\em extended Euclidean algorithm} produces a sequence $(u_i)$
with $u_{-1}:=0$, $u_0:=1$, and $u_i:=Q_iu_{i-1}+u_{i-2}$.
Then for each $i$, $(-1)^ir_i/u_i\equiv r_0 (mod\ N_L)$.
Of these there is necessarily some $i\geq 0$ with $r_i^2+u_i^2$ minimum,
choosing $i$ minimum as well if this is not unique.

Now define the composite map 
\[ \psi^{(0,N_L)}:=
E_{N_L}\circ CRT \circ\left(\prod \{\psi^{(q_l)}: l\in L\}\right) 
\circ\left( \prod\{ \mu_{q_l}: l\in L\}\right)\]
for $\psi^{(q_l)}$ the corresponding inclusion map from $R^{(q_l)}$ to $\ov{R}^{(q_l)}$.

Suppose the variables $\ov{y}_k^{(0)}$ in the integral closure presentation
$\ov{S}^{(0)}$ correspond to the fractions $g_k^{(0)}/\Delta^{(0)}$ 
for $g_k^{(0)},\Delta^{(0)}\in {R}^{(0)}$,
and the Gr\"obner basis elements 
of the ideal of induced relations $\ov{I}^{(0)}$ are denoted by $\ov{b}_k^{(0)}$.
Let the variable $\ov{y}_j^{(q)}$ correspond to $g_j^{(q)}/\Delta^{(q)}$, 
for $g_j^{(q)}:=\mu_q\left(g_j^{(0)}\right)$
and $\Delta^{(q)}:=\mu_q\left(\Delta^{(0)}\right)$. 
If $q$ is a good prime, then these should be variables 
and (a Gr\"obner basis of) relations for the integral closure $\ov{S}^{(q)}$.

Since the object here is to go in the reverse direction by 
reconciling various presentations, $\ov{S}^{(q)}$, 
and reconstructing the presentation $\ov{S}^{(0)}$ from them,
using the Chinese remainder map and the extended Euclidean algorithm map,
consider the candidates for $\ov{S}^{(0)}$, namely $\ov{S}^{(0,N_L)}$ 
with polynomial ring  $\ov{R}^{(0,N_L)}$ having variables $\ov{y}_j^{(0,N_L)}$
corresponding to $g_j^{(0,N_L)}/\Delta^{(0)}$ for
\[ g_j^{(0,N_L)}:=E_{N_L}\circ CRT\circ \prod\{g_j^{(q_l)}:\ l\in L\}\]
and ideal $\ov{I}^{(0,N_L)}$ generated by the {\em finite ordered set} of images
\[\ov{b}_k^{(0,N_L)}:=E_{N_L}\circ CRT \circ \prod\{\ov{b}_k^{(q_l)}:\ l\in L\}.\]

\newpage
\section{Theory}

\begin{lem} If $q=N_1$ is a good prime larger than $\al^2+\be^2$ 
for any coefficient $\al/\be\in\mathbf{Q}$ needed to be reconstructed 
to produce the presentation $\ov{R}^{(0)}/\ov{I}^{(0)}$, 
then $\ov{R}^{(q)}/\ov{I}^{(q)}$ lifts to this presentation.
$[$And the canonical polynomial set $(g_j^{(q)}:0\leq j\leq J)$ necessarily lifts to
a canonical polynomial set $(g_j^{(0,q))}: 0\leq j\leq J)$.$]$
\end{lem}

\begin{pf} If $g_j^{(q)}$ lifts to $g_j^{(0)}$ 
(including $g_0^{(q)}=\Delta^{(q)}$ lifting to $g_0^{(0)}=\Delta^{(0)}$), 
and the relations $\ov{b}_j^{(q)}$ lift to $\ov{b}_j^{(0)}$, 
then $\ov{S}^{(q)}$ lifts to $\ov{S}^{(0)}$.
But if $q>\al^2+\be^2$ then $c\equiv \al/\be(mod\ q)$ lifts to $\al/\be$ 
using the extended Euclidean algorithm as described above.
\end{pf}

\begin{cor} If $(q_l)_{l\in L}$ is a set of distinct good primes and $N_L:=\prod \{q_l:\ l\in L\}$ 
is larger than $\al^2+\be^2$ for any rational coefficient needed to be reconstructed 
to produce the presentation  $\ov{R}^{(0)}/\ov{I}^{(0)}$,  
and $\ov{R}^{(q_l)}/\ov{I}^{(q_l)}$ are compatible in the sense that $LM(g_j^{(q_l)})$ is independent of $l$
and $LM(\ov{b}_k^{(q_l)})$ is independent of $l$, 
then $(\ov{R}^{(q_l)}/\ov{I}^{(q_l)})_{l\in L}$ lifts to this presentation. 
$[$And, again, the canonical polynomial set $(g_j^{(N_L)}:0\leq j\leq J)$ necessarily lifts to
a canonical polynomial set $(g_j^{(0,N_L)}: 0\leq j\leq J)$.$]$
\end{cor}

\begin{pf} Use the Chinese remainder theorem to reconcile these individual presentations, 
and lift the resulting ordered sets $(g_j^{(N_L)})$, $1\leq j\leq J_L$, 
and $(\ov{b}_k^{(N_L)})$, $1\leq k\leq K_L$, 
to ordered sets $(g_j^{(0,N_L)})$, $0\leq j\leq J_L$,  and $(\ov{b}_k^{(0,N_L)})$, $1\leq k\leq K_L$,
and proceed as in the previous proposition.
\end{pf}
Since $\ov{S}^{(0)}$ is not known ahead of time, 
it is not clear how big $N_L$ must be to apply the proposition or corollary above. 
It is therefore better to have a theorem independent of this knowledge.
So the following is a way of knowing that $N_L$ is sufficiently large 
without knowing just how large sufficiently large is.

For $\ov{S}^{(0,N)}$ to be a presentation of the integral closure of $S^{(0)}$, 
it necessarily must be a ring containing $S^{(0)}$ 
and also contained in $\frac{1}{\Delta^{(0)}}S^{(0)}$,
$\ov{S}^{(0)}$ being (isomorphic to) the union of all such. 

\begin{thm} If $\ov{B}^{(0,N)}$ is a Gr\"obner basis for $\ov{I}^{(0,N)}$,
and $\psi^{(N)}(I^{(0)})\subseteq\ov{I}^{(0,N)}$, then 
 $\ov{S}^{(0,N)}=\psi^{(N)}(\ov{S}^{(0)})$. 
\end{thm}

\begin{pf} 
If $\ov{B}^{(0,N)}$ is a Gr\"obner basis for $\ov{I}^{(0,N)}$,
then the quotient ring $\ov{S}^{(0,N)}:=\ov{R}^{(0,N)}/\ov{I}^{(0,N)}$
is a strict affine $P^{(0)}$-algebra.

If $\psi^{(N)}(I^{(0)})\subseteq\ov{I}^{(0,N)}$, then 
\[\psi^{(N)}(S^{(0)})\subseteq\ov{S}^{(0,N)}\subseteq \frac{1}{\Delta^{(0)}}\psi^{(N)}(S^{(0)}).\]
But $\psi^{(N)}(\ov{S}^{(0)})$ is the union of all such rings, so
$\ov{S}^{(0,N)}\subseteq\psi^{(N)}(\ov{S}^{(0)})$. 

If $\ov{S}^{(0,N)}\neq\psi^{(N)}(\ov{S}^{(0)})$, 
consider the monic conductor element $\Delta^{(0)}\in P^{(0)}$ 
mapping to $\Delta^{(q)}\in P^{(q)}$
for all primes $q$ not identified as bad primes 
by the Jacobian computation above.
Were the integral closure of $\psi^{(N)}(S^{(0)})$ computed
as the integral closure of $\ov{S}^{(0,N)}$,
the conductor element $\Delta^{(0,N)}\in P^{(0)}$ computed
would necessarily be a divisor of $\Delta^{(0)}$
since $\ov{S}^{(0,N)}\subset\psi^{(N)}(\ov{S}^{(0)})$.
So $\Delta^{(0,N)}$ would be monic 
with mod $q$ image $\mu_q(\Delta^{(0,N)})\in P^{(q)}$ for any $q|N$. 
But $\mu_q(\ov{S}^{(0,N)})=\ov{S}^{(q)}$ for these $q$.
Since $\ov{S}^{(q)}$ is integrally closed, $\mu_q(\Delta^{(0,N)})=1$.
Hence $\Delta^{(0,N)}$, being a monic element of $P^{(0)}$, 
would be $1$ as well; 
meaning that $\ov{S}^{(0,N)}$ would be integrally closed.
 

\end{pf}

\begin{exmp} Let
\[ S^{(q)}:=\mathbf{F}_q[y_7;x_3]/\langle y_7^3+x_3^7+8y_7x_3\rangle,\]
$($with $\mathbf{F}_0:=\mathbf{Q}$ allowed$)$. 
Then 
\[ \ov{S}^{(0)}=
\mathbf{Q}[z_{11},y_7;x_3]/\langle z_{11}^2+y_7x_3^5+8z_{11},\ z_{11}y_7+x_3^6+8y_7,\ y_7^2-z_{11}x_3\rangle,\]
with the subscripts defining an induced weight function 
and a corresponding induced weight-over-grevlex monomial ordering.
For any prime $q\neq 2$, $\mu_q\left(\ov{S}^{(0)}\right)=\ov{S}^{(q)}$.
So, for any $N$ a product of distinct odd primes, 
$\ov{S}^{(0,N)}=\ov{S}^{(0)}$ if $N> 8^2+1^2$
For a smaller $N$ such as $N=55$, $8(mod\ 55)$ lifts to possibly the wrong fraction, here $1/7$ instead of $8/1$,
and subsequently
\[ y_7^3+x_3^7+8y_7x_3=y_7(y_7^2-z_{11}x_3)+x_3(z_{11}y_7+x_3^6+\frac{1}{7}y_7)+\frac{55}{7}y_7x_3
\notin\ov{I}^{(0,55)}.\]

For $q=2$,
\[\ov{S}^{(2)}=\mathbf{F}_2[w_1;x_3]/\langle w_1^3+x_3\rangle,\]
with $y_7=w_1x_3^2$ and $z_{11}=w_1^2x_3^3$. 
Clearly this is larger than expected, so $\ov{S}^{(2)}\supset \ov{S}^{(0)}$.
\end{exmp}

\begin{exmp}
The generic example $\mathbf{F}_q[y_1;x_1]/\langle y_1^3+a_qy_1x_1+b_qx_1^5\rangle$
has $y_1^2/x_1$ in its integral closure;
so its integral closure has a presentation as
\[\mathbf{F}_q[\ov{y}_2,\ov{y}_1;x_1]/
\langle \ov{y}_2^2+a_q\ov{y}_2+b_q\ov{y}_1x_1^3,\ \ov{y}_2\ov{y}_1+a_q\ov{y}_1+b_qx_1^4,\ \ov{y}_1^2-\ov{y}_2x_1\rangle.\]

If $\mathbf{F}_0=\mathbf{Q}$, and $a_0:=1/3$ and $b_0:=8/7$, then  
the image in characteristic $q$ is not defined for $q=3,7$, and is not an affine domain for $q=2$.
For $\mathbf{F}_5=\mathbf{Z}_5$, $a_5=2$ and $b_5=-1$ would lift to $a_0=2$ and $b_0=-1$, 
giving a presentation of the wrong integral closure $($one with the right form but these wrong coefficients$)$.
Using $\mathbf{F}_{11}=\mathbf{Z}_{11}$ as well would give $a_{11}=4$ and $b_{11}=-2$ 
reconciled to give $a_{55}=-18$ and $b_{55}=9$,
and lifted to $a_0=1/3$ and $b_0=-1/6$,
again giving a presentation of the wrong integral closure.
Using in addition $\mathbf{F}_{13}=\mathbf{Z}_{13}$ would produce $a_{13}=-3$ and $b_{13}=3$ 
reconciled to give $a_{715}=-238$ and $b_{715}=-101$,
lifted to the correct $a_0=1/3$ and $b_0=8/7$.
\end{exmp}

\begin{exmp}
$\mathbf{Q}[y,x]/\langle y^2-3/2x^3+24/7x^2-96/49x\rangle$ doesn't need $N> 96^2+49^2$ to work, but only
$N>8^2+7^2$, since the only things needed to be computed are $g_0:=x-8/7$, $g_1:=y$, and $\ov{b}_1:=\ov{y}^2-3/2x$
$($and the inclusion map image $\psi(y):=\ov{y}(x-8/7)$$)$.

The details for this example are as follows:
The primes $q=2,7$ are bad because they divide denominators of fractions defining the problem.
The image for $q=3$ is not even a reduced ring, so probably should be avoided as well.
$\delta^{(0)}=7x-8$, and $q=7$ is already to be avoided.

For $q=5$, $g_0^{(5)}=\Delta^{(5)}=x+1$, $g_1^{(5)}=y$, $\ov{b}^{(5)}=\ov{y}^2+x$, $\psi^{(5)}(y)=\ov{y}(x+1)$.
This lifts to give $g_0^{(0,5)}=x+1$, $g_1^{(0,5)}=y$, $\ov{b}_1^{(0,5)}=\ov{y}^2+x$, $\psi^{(5)}(y)=\ov{y}(x+1)$.
Then the defining relation above reduces to $-x(x+1)^2 -3/2x^3+24/7x^2-96/49x\neq 0$.

For $q=11$, $g_0^{(11)}=\Delta^{(11)}=x+2$, $g_1^{(11)}=y$, $\ov{b}^{(11)}=\ov{y}^2-4x$, $\psi^{(11)}(y)=\ov{y}(x+2)$.
This reconciles with the previous to get $g_0^{(55)}=x-9$, $g_1^{(55)}=y$, $\ov{b}^{(55)}=\ov{y}^2+26x$, $\psi^{(55)}(y)=\ov{y}(x-9)$.
This lifts to give $g_0^{(0,55)}=x+1/6$, $g_1^{(0,55)}=y$, $\ov{b}_1^{(0,55)}=\ov{y}^2-3/2x$, $\psi^{(55)}=\ov{y}(x+1/6)$.
Then the defining relation above reduces to $3/2x(x+1/6)^2 -3/2x^3+24/7x^2-96/49x\neq 0$.

For $q=13$, $g_0^{(13)}=\Delta^{(13)}=x-3$, $g_1^{(13)}=y$, $\ov{b}^{(13)}=\ov{y}^2+5x$, $\psi^{(13)}(y)=\ov{y}(x-3)$.
This reconciles with the previous to get $g_0^{(715)}=x+101$, $g_1^{(715)}=y$, $\ov{b}^{(715)}=\ov{y}^2+356x$, $\psi^{(715)}(y)=\ov{y}(x+101)$.
This lifts to give $g_0^{(0,715)}=x-8/7$, $g_1^{(0,715)}=y$, $\ov{b}_1^{(0,715)}=\ov{y}^2-3/2x$, $\psi^{(715)}=\ov{y}(x-8/7)$.
Then the defining relation above reduces to $3/2x(x-8/7)^2 -3/2x^3+24/7x^2-96/49x=0$.

The presentation found $($but not minimized$)$ is then
\[\mathbf{Q}[\ov{y};x]/\langle \ov{y}^2-3/2x\rangle\]
with inclusion map defined by $\psi(y)=\ov{y}(x-8/7)$.
$[$ The minimized presentation here would have been just the polynomial ring $\mathbf{Q}[\ov{y}]$
with $x=2/3\ov{y}^2$ and $y=\ov{y}(2/3\ov{y}^2-8/7)$ both unnecessary except for defining the inclusion.$]$
\end{exmp}

\newpage

\section{Appendix}
It is envisioned that the code and the relevant examples relative to this paper
on the website {\tt http://www.dms.auburn.edu/\~{}leonada}.
will be updated as various packages change for the better.
The code for the Qth-power algorithm in positive characteristic
and the extra code to extend it to char $0$ for this paper
are both written in {\sc Magma} and in {\sc Macaulay2}
and are available from the author.

But, as mentioned in the Overview section above, 
the complete version of the example mentioned there,
is done here by the various methods mentioned.
  
$P^{(0)}:=\mathbf{Q}[x]$,
$R^{(0)}:=\mathbf{Q}[y;x]$, 
\[B^{(0)}:=\{(y^2-3/4y-15/17x)^3-9yx^4(y^2-3/4y-15/17x)-27x^{11}\}\]

The Qth-power algorithm implementation produces fractions with numerators
\begin{verbatim}
    p_5:=x^5, 
    p_4:=y^2*x^3-(3/4)*y*x^3-(15/17)*x^4, 
    p_3:=y*x^5,
    p_2:=y^4*x-(3/2)*y^3*x-(30/17)*y^2*x^2+(9/16)*y^2*x+(45/34)*y*x^2
       +(225/289)*x^3, 
    p_1:=y^3*x^3-(15/17)*y*x^4-(9/16)*y*x^3-(45/68)*x^4,
    p_0:=y^5-(9/4)*y^4-(30/17)*y^3*x+(27/16)*y^3+(45/17)*y^2*x
       +(225/289)*y*x^2-(27/64)*y^2-(135/136)*y*x-(675/1156)*x^2
\end{verbatim}
$p_5$ being the common denominator, a conductor element lying in $P$, though $\Delta^{(0)}=x^9$ is the one computed
directly from the Jacobian.
The implementation also produces a Gr\"obner basis $\ov{B}$ for the presentation:
\begin{verbatim}
     p_0^2-(135/17)*p_0+(81/4)*p_1*p_5^3-27*p_2*p_5^5-81*p_2*p_5^2
        -243*p_3*p_5^5-(405/17)*p_4*p_5^4-(243/8)*p_4*p_5^3
        -(1215/17)*p_4*p_5+(729/4)*p_5^5,
     p_0*p_1-9*p_0*p_5^2-(135/17)*p_1-(27/2)*p_2*p_5
        -(81/4)*p_3*p_5^4-27*p_4*p_5^6-(405/17)*p_5^5
        +(243/16)*p_5^4,
     p_0*p_2-27*p_1*p_5^4-81*p_1*p_5-(135/17)*p_2+(81/2)*p_4*p_5^4
        +(243/4)*p_4*p_5-243*p_5^6, 
     p_0*p_3-9*p_1*p_5-(15/17)*p_2+(27/4)*p_4*p_5-27*p_5^6,
     p_0*p_4-9*p_2*p_5-27*p_3*p_5^4-(135/17)*p_4+(81/4)*p_5^4,
     p_1^2-(9/4)*p_0*p_5-9*p_1*p_5^2-(15/17)*p_2*p_5-(9/4)*p_2
        +(27/4)*p_4*p_5^2-27*p_5^7, 
     p_1*p_2-(27/2)*p_1-9*p_2*p_5^2-27*p_3*p_5^5-(135/17)*p_4*p_5
        +(81/8)*p_4-(81/4)*p_5^5,
     p_1*p_3-(3/2)*p_1-p_2*p_5^2-(15/17)*p_4*p_5+(9/8)*p_4,
     p_1*p_4-p_0*p_5-(3/2)*p_2, 
     p_2^2-9*p_0*p_5-(27/4)*p_2-27*p_4*p_5^5,
     p_2*p_3-p_0*p_5-(3/4)*p_2, p_2*p_4-9*p_1+(27/4)*p_4-27*p_5^5,
     p_3^2-(3/4)*p_3-p_4*p_5^2-(15/17)*p_5, 
     p_3*p_4-p_1+(3/4)*p_4, 
     p_4^2-p_2
\end{verbatim}
with induced weights $wt(p_0)=25$, $wt(p_1)=21$, $wt(p_2)=20$, $wt(p_3)=11$, $wt(p_4)=10$, and $wt(p_5)=6$.
[Note that $(p_0,p_1,p_2,p_3,p_4,p_5)$ here correspond 
to $(\ov{y}_4,\ov{y}_3,\ov{y}_2,\ov{y}_1,x_1)$ in the notation of this paper;
but the notation $(f_{25},f_{21},f_{20},f_{11},f_{10},f_6)$ would be better than either, 
given that the subcripts then correspond to the weights.]

Using {\sc Macaulay2}'s {\em integralClosure} function, \cite{M2}, 
an implementation of de Jong's algorithm, \cite{Jong},
the output ideal is generated by:
\begin{verbatim}
  314432y6-8489664x11-2829888y3x4-707472y5-832320y4x+2122416y2x4
    +530604y4+2496960yx5+1248480y3x+734400y2x2-132651y3-468180y2x
    -550800yx2-216000x3,
  w_(3,0)x2-68y2+51y+60x, 
  4624w_(3,0)y4-6936w_(3,0)y3-8160w_(3,0)y2x+2601w_(3,0)y2
     +6120w_(3,0)yx-8489664x9-2829888y3x2+2122416y2x2
     +2496960yx3+244800y2-183600y-216000x,
  68w_(3,0)^2y2-51w_(3,0)^2y-60w_(3,0)^2x-8489664x7-2829888y3
     +2122416y2+2496960yx w_(3,0)^3-41616w_(3,0)y-8489664x5,
  w_(4,0)x-4w_(3,0)^2y+3w_(3,0)^2,
  17w_(4,0)y-60w_(3,0)^2-41616w_(3,0)yx-8489664x6,
  w_(4,0)w_(3,0)-2448w_(3,0)^2x-146880w_(3,0)-33958656yx4
     +25468992x4,
  w_(4,0)^2-146880w_(4,0)-407503872w_(3,0)y2+305627904w_(3,0)y
     -9236754432y4x-83130789888yx5+20782697472y3x+8150077440y2x2
     +62348092416x5-15587023104y2x-12225116160yx2
     +3896755776yx+4584418560x2|
\end{verbatim}

This is a presentation relative to $\mathbf{Z}$, as attested to by the leading coefficients. 
One can use $w_{3,0}=68b$ and $w_{4,0}=17*68^2a$ to clean this up a bit, 
but it will still be an affine $S$-algebra presentation with a default block ordering, 
grevlex on the new variables, forced to have the input over $\mathbf{Z}$ as an explicit subring.
The fact that it essentially found a single common denominator conductor element in $P$ is uncharacteristic.

Trying {\sc Singular}'s {\em normal} function, \cite{GPS}\cite{GLS}, 
also an implementation of de Jong's algorithm
gives numerators:
\begin{verbatim}
   68y2x3-51yx3-60x4,
   4624y4x-6936y3x-8160y2x2+2601y2x+6120yx2+3600x3,
   18496y5-41616y4-32640y3x+31212y3+48960y2x+14400yx2-7803y2
      -18360yx-10800x2,
   x5
\end{verbatim}
The relations are:
\begin{verbatim}
s[ 1]=314432*y^6-8489664*x^11-2829888*y^3*x^4-707472*y^5
     -832320*y^4*x+2122416*y^2*x^4+530604*y^4+2496960*y*x^5
     +1248480*y^3*x+734400*y^2*x^2-132651*y^3-468180*y^2*x
     -550800*y*x^2-216000*x^3
s[ 2]=68*T(3)*y^2*x-51*T(3)*y*x-60*T(3)*x^2-33958656*y*x^7
     -11319552*y^4+25468992*x^7+16979328*y^3+9987840*y^2*x
     -6367248*y^2-7490880*y*x
s[ 3]=3468*T(3)*y*x^3-340*T(3)*y^2+255*T(3)*y+300*T(3)*x
     -1731891456*x^9+169793280*y*x^6-577297152*y^3*x^2
     -127344960*x^6+432972864*y^2*x^2+509379840*y*x^3
     +49939200*y^2-37454400*y-44064000*x
s[ 4]=T(3)*x^5-18496*y^5+41616*y^4+32640*y^3*x-31212*y^3
     -48960*y^2*x-14400*y*x^2+7803*y^2+18360*y*x+10800*x^2
s[ 5]=272*T(3)*y^3-2448*T(3)*x^4-408*T(3)*y^2-240*T(3)*y*x
     +153*T(3)*y+180*T(3)*x-135834624*y^2*x^6+203751936*y*x^6
     -76406976*x^6-39951360*y^3+59927040*y^2+35251200*y*x
     -22472640*y-26438400*x
s[ 6]=60*T(2)*x-17*T(3)*y*x+8489664*x^7+2829888*y^3
     -2122416*y^2-2496960*y*x
s[ 7]=4*T(2)*y-3*T(2)-T(3)*x
s[ 8]=T(1)*x^2-68*y^2+51*y+60*x
s[ 9]=41616*T(1)*y*x+60*T(2)-17*T(3)*y+8489664*x^6
s[10]=4080*T(1)*y^2-3060*T(1)*y-3600*T(1)*x-17*T(3)*y*x^2
     +8489664*x^8+2829888*y^3*x-2122416*y^2*x-2496960*y*x^2
s[11]=5*T(3)^2-1797811200*T(1)*x-8489664*T(3)*y*x^2
     -734400*T(3)-46183772160*y^4*x+4239670284288*x^8
     -415653949440*y*x^5+1517136915456*y^3*x+40750387200*y^2*x^2
     +311740462080*x^5-1137852686592*y^2*x-1308087429120*y*x^2
     +19483778880*y*x+22922092800*x^2
s[12]=T(2)*T(3)-41616*T(3)*y-2309188608*y^3*x^2
     +3463782912*y^2*x^2+2037519360*y*x^3-1298918592*y*x^2
     -1528139520*x^3
s[13]=5*T(1)*T(3)-734400*T(1)-3468*T(3)*y*x+1731891456*x^7
      -169793280*y*x^4
     +577297152*y^3+127344960*x^4-432972864*y^2-509379840*y*x
s[14]=T(2)^2-31212*T(2)-10404*T(3)*x-577297152*y^2*x^3
     +432972864*y*x^3+509379840*x^4
s[15]=T(1)*T(2)-41616*T(1)*y-8489664*x^5
s[16]=T(1)^2-T(2)
\end{verbatim}

Again, this is written relative to $\mathbf{Z}$
and can be cleaned up a bit by using $T(1)=68*c$, $T(2)=68^2*b$, $T(3)=17*68^2*a$
This is at least a strict affine $S$-algebra presentation, but again suffers
from being relative to $S$ and having no hint of the induced monomial ordering.
The fact that it found a conductor element in $P$ is uncharacterisitc. 

In {\sc Magma}, \cite{Magma},  the {\em Normalisation} function, 
a third implementation of de Jong's algorithm, gives a basis:

\begin{verbatim}
[$.1^4-4913/3375*$.1*$.3^6+4913/1500*$.1*$.3^5-4913/2000*$.1*$.3^4
   +4913/8000*$.1*$.3^3+17/405*$.2^4-4913/10125*$.2^3*$.3^3
   +4913/13500*$.2^3*$.3^2+289/6075*$.2^2*$.3*$.4-17/45*$.2^2*$.3
   +4913/91125*$.3^2*$.4^2-578/675*$.3^2*$.4,
 $.1^3*$.3-3/4*$.1^3-289/225*$.1*$.3^5+289/100*$.1*$.3^4
   -867/400*$.1*$.3^3+867/1600*$.1*$.3^2-289/675*$.2^3*$.3^2
   +289/900*$.2^3*$.3+17/405*$.2^2*$.4+289/6075*$.3*$.4^2
   -34/45*$.3*$.4,
 $.1^2*$.2+15/17*$.1-$.3^2+3/4*$.3,
 $.1^2*$.3^2-3/2*$.1^2*$.3+9/16*$.1^2-17/15*$.1*$.3^4
   +51/20*$.1*$.3^3-153/80*$.1*$.3^2+153/320*$.1*$.3
   -17/45*$.2^3*$.3+17/60*$.2^3+17/405*$.4^2-1/3*$.4,
 $.1*$.2^2+4913/1125*$.1*$.3^7-4913/375*$.1*$.3^6
   +14739/1000*$.1*$.3^5-14739/2000*$.1*$.3^4
   +44217/32000*$.1*$.3^3-17/135*$.2^4*$.3+17/180*$.2^4
   +4913/3375*$.2^3*$.3^4-4913/2250*$.2^3*$.3^3
   +4913/6000*$.2^3*$.3^2-289/2025*$.2^2*$.3^2*$.4
   +17/15*$.2^2*$.3^2+289/2700*$.2^2*$.3*$.4-17/20*$.2^2*$.3
   -1/9*$.2*$.4+15/17*$.2-4913/30375*$.3^3*$.4^2
   +578/225*$.3^3*$.4+4913/40500*$.3^2*$.4^2-289/150*$.3^2*$.4,
 $.1*$.2*$.3^7-3*$.1*$.2*$.3^6+27/8*$.1*$.2*$.3^5
   -27/16*$.1*$.2*$.3^4+81/256*$.1*$.2*$.3^3-25/867*$.2^5*$.3
   +25/1156*$.2^5+1/3*$.2^4*$.3^4-1/2*$.2^4*$.3^3
   +3/16*$.2^4*$.3^2-5/153*$.2^3*$.3^2*$.4+75/289*$.2^3*$.3^2
   +5/204*$.2^3*$.3*$.4-225/1156*$.2^3*$.3-125/4913*$.2^2*$.4
   -1/27*$.2*$.3^3*$.4^2+10/17*$.2*$.3^3*$.4+1/36*$.2*$.3^2*$.4^2
   -15/34*$.2*$.3^2*$.4+1125/4913*$.3*$.4,
 $.1*$.3^14-6*$.1*$.3^13+63/4*$.1*$.3^12-189/8*$.1*$.3^11
   +2835/128*$.1*$.3^10-1701/128*$.1*$.3^9+5103/1024*$.1*$.3^8 
   -2187/2048*$.1*$.3^7+6561/65536*$.1*$.3^6
   +9375/1419857*$.2^6*$.3-28125/5679428*$.2^6
   -375/4913*$.2^5*$.3^4+1125/9826*$.2^5*$.3^3 
   -3375/78608*$.2^5*$.3^2-25/867*$.2^4*$.3^8
   +125/1156*$.2^4*$.3^7-375/2312*$.2^4*$.3^6
   +1125/9248*$.2^4*$.3^5-3375/73984*$.2^4*$.3^4 
   +2025/295936*$.2^4*$.3^3+625/83521*$.2^4*$.3^2*$.4
   -84375/1419857*$.2^4*$.3^2-1875/334084*$.2^4*$.3*$.4 
   +253125/5679428*$.2^4*$.3+1/3*$.2^3*$.3^11-3/2*$.2^3*$.3^10 
   +45/16*$.2^3*$.3^9-45/16*$.2^3*$.3^8+405/256*$.2^3*$.3^7
   -243/512*$.2^3*$.3^6+243/4096*$.2^3*$.3^5
   +140625/24137569*$.2^3*$.4-5/153*$.2^2*$.3^9*$.4
   +75/289*$.2^2*$.3^9+25/204*$.2^2*$.3^8*$.4
   -1125/1156*$.2^2*$.3^8-25/136*$.2^2*$.3^7*$.4
   +3375/2312*$.2^2*$.3^7+75/544*$.2^2*$.3^6*$.4
   -10125/9248*$.2^2*$.3^6-225/4352*$.2^2*$.3^5*$.4
   +30375/73984*$.2^2*$.3^5+135/17408*$.2^2*$.3^4*$.4
   -18225/295936*$.2^2*$.3^4+125/14739*$.2^2*$.3^3*$.4^2
   -11250/83521*$.2^2*$.3^3*$.4-125/19652*$.2^2*$.3^2*$.4^2
   +16875/167042*$.2^2*$.3^2*$.4-125/4913*$.2*$.3^7*$.4
   +16875/83521*$.2*$.3^7+375/4913*$.2*$.3^6*$.4
   -50625/83521*$.2*$.3^6-3375/39304*$.2*$.3^5*$.4
   +455625/668168*$.2*$.3^5+3375/78608*$.2*$.3^4*$.4
   -455625/1336336*$.2*$.3^4-10125/1257728*$.2*$.3^3*$.4
   +1366875/21381376*$.2*$.3^3-1265625/24137569*$.2*$.3*$.4
   -1/27*$.3^10*$.4^2+10/17*$.3^10*$.4+5/36*$.3^9*$.4^2
   -75/34*$.3^9*$.4-5/24*$.3^8*$.4^2+225/68*$.3^8*$.4
   +5/32*$.3^7*$.4^2-675/272*$.3^7*$.4-15/256*$.3^6*$.4^2
   +2025/2176*$.3^6*$.4+9/1024*$.3^5*$.4^2-1215/8704*$.3^5*$.4,
 $.1*$.4 - $.2^2*$.3 + 3/4*$.2^2]
\end{verbatim}

But the computed Gr\"obner basis for this is way to big to be reproduced here,
since {\em Normalisation} almost always chooses a default {\em Lex} monomial ordering the reverse
of that which will give a readable {\em Lex} answer. 
As case-in-point, one can get a readable {\em Lex} monomial order answer 
by reversing the variables before computing a Gr\"obner basis:

\begin{verbatim}
 $.1^2-135/17*$.1-$.3^5+9/4*$.3^3*$.4^3+27/4*$.3^3-405/17*$.3*$.4^4
    -243/16*$.3*$.4^3-243/4*$.4^8,
 $.1*$.2-$.3^3*$.4-15/17*$.3^2,
 $.1*$.3-27*$.2*$.4^4-9*$.3^2*$.4-135/17*$.3+81/4*$.4^4,
 $.1*$.4-1/9*$.3^4+3/4*$.3^2+3*$.3*$.4^5,
 $.2^2-3/4*$.2-$.3*$.4^2-15/17*$.4,
 $.2*$.3-1/9*$.3^3+3*$.4^5,
 $.2*$.4^5-1/243*$.3^5+1/36*$.3^3+1/9*$.3^2*$.4^5+1/3*$.3^2*$.4^2
    +5/17*$.3*$.4-3/4*$.4^5,
 $.3^6-27/4*$.3^4-54*$.3^3*$.4^5-81*$.3^3*$.4^2-1215/17*$.3^2*$.4 
    +729/4*$.3*$.4^5+729*$.4^10
\end{verbatim} 

Since this implementation does not force a presentation relative to S,
it occassionally gives a decent minimized presentation.
But there is, again, no hint that there is a natural induced monomial ordering.

Since there is only one free variable in this example, {\sc Magma}'s
{\em IntegralClosure} gives an answer

\begin{verbatim} 
 B[1]= 1
 B[2]= Y
 B[3]= 1/X^2*Y^2-3/4/X^2*Y-15/17/X
 B[4]= 1/X^2*Y^3+(-15/17*X-9/16)/X^2*Y-45/68/X
 B[5]= 1/X^4*Y^4-3/2/X^4*Y^3+(-30/17*X+9/16)/X^4*Y^2+45/34/X^3*Y
       +225/289/X^2
 B[6]= 1/X^5*Y^5-9/4/X^5*Y^4+(-30/17*X+27/16)/X^5*Y^3
       +(45/17*X-27/64)/X^5*Y^2+(225/289*X-135/136)/X^4*Y
       -675/1156/X^3
\end{verbatim}

At least this necessarily gives a $P$-module basis 
and an answer over $\mathbf{Q}$ instead of $\mathbf{Z}$. 
But there is obviously no way to give weights,
and the presentation is only implicit.

\end{document}